\documentclass[a4paper]{easychair}

\usepackage{stmaryrd} 
\usepackage{amssymb}
\newcommand\N{\mathbb{N}}
\usepackage{amsmath}
\usepackage{amsfonts}
\usepackage[noxy]{virginialake}

\usepackage{doc}

\usepackage[all]{xy}

\newcommand{\ee}{\varepsilon}

\newcommand{\C}{\textbf{C}}
\newcommand{\Pos}{\textbf{Pos}}
\newcommand{\CTX}{\textbf{Ctx}}
\newcommand{\LT}{\textsf{LT}}

\newtheorem{definition}{Definition}
\newtheorem{proposition}{Proposition}

\title{First-order Intuitionistic Multiplicative Linear Logic with Hyperdoctrines and Hilbert's Epsilon\thanks{AILA 2026, XXIX Incontro di Logica, Rome, 15--18 September 2026, Booklet of Abstracts, pp.~30--33.}}

\author{
   Loïc Allègre\inst{1}
\and 
Fabio Pasquali\inst{2}
\and
    Christian Retoré\inst{1}
}

\institute{LIRMM - University of Montpellier/CNRS,
  Montpellier, France\\
  \email{loic.allegre@lirmm.fr\quad christian.retore@lirmm.fr}
\and 
  Università degli Studi di Milano,
   Milano, Italy\\
   \email{fabio.pasquali@unimi.it}\\
 }
 
\authorrunning{Allègre, Pasquali, Retoré}

\titlerunning{Hyperdoctrines for $\epsilon$ First-order Intuitionistic Linear Logic}

\begin{document}

\maketitle

\begin{abstract}
We propose an IMELL doctrine whose base category is the full subcategory of Set consisting of
finite powers of a given set $U$,  using a particular monoidal poset $M^{\omega}$ with a non-trivial comonad. This  doctrine is shown to be an $\epsilon$-doctrine, and thus provides a categorical
semantics for  IMELL $\epsilon$-calculus, a first order intuitionistic linear logic.  
\end{abstract}

% The table of contents below is added for your convenience. Please do not use
% the table of contents if you are preparing your paper for publication in the
% EPiC Series or Kalpa Publications series

%\setcounter{tocdepth}{2}
%{\small
%\tableofcontents}

%\paragraph{To mention}
%
%Processing in EasyChair - number of pages.
%
%Examples of how EasyChair processes papers. Caveats (replacement of EC
%class, errors).

%------------------------------------------------------------------------------
\paragraph{Introduction}
\label{sect:introduction}
Girard's Linear Logic \cite{Gir87} has established itself as an exceptionally precise framework for the fine-grained analysis of classical logic and, most notably, intuitionistic logic. Beyond its logical significance, it has also proved highly fruitful as a logic of computation. In particular, it has played a central role in the study of cut-elimination and $\beta$-reduction, most famously through the demonstration that Lamping’s algorithm realises maximal sharing,\cite{GonthierAbadiLevy1992OptimalReduction} as well as through the development of Geometry of Interaction \cite{Girard2011blindspot} 
and has provided powerful tools for the modelling of concurrent processes, as recently exemplified in the work by Manara \cite{Manara2025LinearLogicPiCalculus}.

From a logical perspective, however, Linear Logic has largely left aside the systematic study of first order quantification — as in typed lambda calculus,  second order quantification has been more successful 
\cite{Girard2011blindspot}.
Standard first-order quantifiers can of course be incorporated into Linear Logic via the usual rules, or even in proof nets as studied in the works of Heijltjes, Hughes and Straßburger \cite{heijltjes_proof_2019}. However, while this extension is technically well behaved, it remains conceptually unsatisfactory. One would expect Linear Logic to articulate a genuinely linear notion of quantification, one that is not merely inherited from classical or intuitionistic logic, but instead reflects the distinctive resource-sensitive character of the theory itself. Indeed, usual quantifiers exhibit an irreducibly additive dimension, even within an otherwise purely multiplicative fragment, a phenomenon already emphasised, for instance, in the work of Fleury-Donnadieu and Quatrini on first-order Ludics \cite{FleuryQuatrini2004FirstOrderLudics}. 

One possible line of approach consists in treating individuals as propositions—understood as individual concepts—while propositions themselves are viewed as classes or terms, and in exploiting second-order quantification accordingly. This perspective is closely related to Girard’s proposal of a transcendental syntax, in which quantification is internalised via higher-order structure rather than expressed through explicit additive operators \cite{Girard2016TranscendentalSyntaxI}.

A more down-to-earth alternative is to turn instead to Hilbert’s subnector,  $\epsilon$ \cite{HBvol2,ChatzikyriakidisPasqualiRetore2017Introduction}, with $\phi[\epsilon_x\, \phi[x]]$ 
expressing $\exists x\ \phi[x]$. Since an expression of the form $\epsilon_x\, \phi[x]$ is a term that effects quantification without appealing to additive connectives, it offers a promising route for reconciling quantification with the multiplicative discipline of Linear Logic. At the same time, we aim to retain truth-valued semantics, so as to preserve a meaningful notion of first-order model. The present work, which remains in progress, is informed by categorical logic for intuitionistic systems, and in particular by the theory of hyperdoctrines \cite{Lawvere1969Adjointness,Pitts1992CategoricalLogic},
which provides a natural semantic environment for such investigations, 
pursuing the work of \cite{DagninoPasquali2025Lipschitz},  on doctrines for linear logic 
and of \cite{Pasquali2017EpsilonDoctrines,Pasquali2016EpsilonIntuitionistic,MaiettiRP} on epsilon in doctrines.

\paragraph{First-order IMELL with Epsilon} 

The language of Intuitionistic Multiplicative Exponential Linear Logic (IMELL) + $\epsilon$ is the usual language of IMELL, extended with the operator $\epsilon$ forming terms $t$ from formulae $\mathcal{F}$:
$$
        t := x \mid c \mid f(t,...,t)\mid \epsilon_x \mathcal{F}
        $$
        $$
        \mathcal{F} := \textbf{1} \mid \bot \mid P(t,...,t) \mid \mathcal{F} \otimes \mathcal{F} \mid \mathcal{F} \multimap \mathcal{F} \mid \ !\mathcal{F}
        $$

The meaning of $\epsilon_x \phi$ is, informally, "some $x$ for which $\phi$ holds, if there are any". This is formally expressed by the corresponding sequent calculus rules, which are the usual rules of IMELL (see \cite{Gir87}), plus the $\epsilon$-rules:
$$
\vlinf{\ee_L}{\text{\scriptsize $x\notin free(\Gamma)$ }}{\Gamma,  A(\ee_x A) \vdash C}{\Gamma, A(x) \vdash C} 
\quad\quad\quad
\vlinf{\ee_R}{}{\Gamma\vdash A(\ee_x A) }{\Gamma \vdash A(t)}
$$

\paragraph{First-order Hyperdoctrines with Epsilon}

Since usual (classical) models of the $\epsilon$-calculus are not adequate for Linear Logic, we instead propose categorical semantics using hyperdoctrines, in the spirit of \cite{Pasquali2016EpsilonIntuitionistic, Pasquali2017EpsilonDoctrines}.

\begin{definition}[\textbf{Hyperdoctrine}]
    A \emph{hyperdoctrine}, or simply doctrine, is a functor $P:\C^{op}\to\Pos$, where $\C$ is a category with finite products (called the base category of $P$), and where $\Pos$ denotes the category of posets and monotone functions.
\end{definition}

The category $\C$ is to be seen as a category of contexts, while the poset $P(A)$ is understood as the collection of formulas over the context $A$, ordered by entailment. Typically, one requires that each $P(A)$ be equipped with certain operations, depending on the logic under consideration. In the case of IMELL, the relevant doctrines are those for which each $P(A)$ is an ordered monoid which, when viewed as a category, is monoidal closed and equipped with a comonad $!_A:P(A)\to P(A)$ satisfying $!(a)\leq !(a) \otimes !(a)$ and $!(a)\leq 1$ \cite{DagninoPasquali2025Lipschitz}.

\begin{definition}[\textbf{Syntactic Doctrine}]
    Let $T$ be an IMELL-theory over a language $L$. The \emph{syntactic doctrine of $T$} is the doctrine $\LT_T:\CTX_L^{op}\to\Pos$, where $\CTX_L^{op}$ is the category of $L$-contexts, in which objects are finite lists of distinct variables, morphisms from $(x_1,\ldots,x_n)$ to $(y_1,\ldots,y_k)$ are lists of $L$-terms $(t_1,\ldots,t_k)$ whose free variables are $(x_1,\ldots,x_n)$, and for all $(x_1,\ldots,x_n)$, the poset $\LT_T(x_1,\ldots,x_n)$ is the Lindenbaum--Tarski algebra of $T$-equiprovable formulas of $L$ whose free variables occur in $(x_1,\ldots,x_n)$.
\end{definition}
Binary products in $\CTX_L^{op}$ are given by context concatenation, and the terminal object is the empty context.  The action of $\LT_T$ on the morphisms of $\CTX_L$ is defined by simultaneous substitution, as depicted in the following diagram.

\centerline{
$
\xymatrix{
(x_1,...,x_n)\ar[d]_-{(t_1,...,t_k)}&\LT_T(x_1,...,x_n)&\phi[t_1/y_1][t_2/y_2]...[t_k/y_k]\\
(y_1,...,y_k)&\LT_T(y_1,...,y_k)\ar[u]_-{\LT_T(t_1,...,t_k)}&\phi\ar@{|->}[u]
}
$} 

The monoidal operations and the comonad on each $\LT_T(x_1,...,x_n)$ are given by the corresponding operations of linear logic.

\begin{definition}[\textbf{$\epsilon$-doctrine}]
    An \emph{$\epsilon$-doctrine} is a doctrine $P:\C^{op}\to\Pos$ such that, for every pair of objects $X,Y$ of $\C$ and every element $\phi$ of $P(X\times Y)$, there exists an arrow $\epsilon_\phi:X\to Y$ satisfying $\exists y\ \phi
\leq P(\langle id_X,\epsilon_\phi\rangle)(\phi)$.
\end{definition}

Here $\exists y\ \phi$ is a simpler notation for $\Sigma(\pi_{X\times Y\mapsto X})(\phi)$ that is the left adjunct to  the image by $P$ of the projection $X\times Y\mapsto X$ in the base category. Accordingly, the class of doctrines associated with the IMELL calculus equipped with the $\epsilon$ subnector is precisely the class of IMELL $\epsilon$-doctrines.  

\paragraph{Categorical Semantics for IMELL + $\epsilon$}

We propose to build for any non-empty set $U$ an $\epsilon$-doctrine $P_U$ whose base category is the full subcategory of $\textbf{Set}$ consisting of finite powers of $U$. We construct the adequate IMELL structure of the posets $P_U(X)$ in the following way.

\begin{definition}
    Let $M$ = $\{0\}\cup \{{\frac{1}{2^n}}|n\in\N\}\subset{\mathbb{Q}}$. We define a \emph{MLL phase structure} as the ordered monoid $(M, \cdot, \leq)$ where $\cdot,\leq$ are the usual multiplication and order on $\mathbb{Q}$.
\end{definition}
We define the operation $\multimap$ on $M$ as $u\multimap v=\min (1,v/u)$ when $u\neq 0$ while $0\multimap v=1$. MLL formulae are interpreted as expected: $\llbracket \textbf{1} \llbracket = 1$, $\llbracket \bot \rrbracket = 0$, $\llbracket A\otimes B \rrbracket=\llbracket A\rrbracket \cdot \llbracket B \rrbracket$ and $\llbracket A\multimap B \rrbracket=\llbracket B\rrbracket \multimap \llbracket A \rrbracket$. 
It is easily checked that $(a\otimes b)\leq c$ iff $a\leq (b \multimap c)$. 

The exponential $!$ is to be interpreted by a comonad, with: 
$\llbracket !A \rrbracket \leq  \llbracket 1 \rrbracket$ (this one holds for any possible value of $\llbracket !A \rrbracket$), 
$\llbracket !A \rrbracket \leq  \llbracket A \rrbracket$, 
$\llbracket !A \rrbracket \leq  \llbracket !A\otimes !A \rrbracket= \llbracket !A \rrbracket \cdot \llbracket !A \rrbracket$. 

In the structure $M$, this leaves only two possible values for $\llbracket !A \rrbracket$: $0$ or $1$, and consequently the interpretation of $!$-formulas is defined by $\llbracket !A \rrbracket= 1 \mbox{ iff } \llbracket A \rrbracket=1$ and otherwise $\llbracket !A \rrbracket= 0$. 

Since this interpretation of ! would be too trivial, we use as phase structure $\omega$ copies of a MLL structure.

\begin{definition}
    We define an \emph{IMELL phase structure} as $(M^\omega, \cdot^\omega, \leq^\omega, !)$, where $M$ is a MLL structure, and the operations on $M^\omega$ are defined pointwise:
    
    \begin{tabular}[t]{l} 
$(u_n)_{n\in\N}\cdot^\omega (v_n)_{n\in\N}=(u_n \cdot v_n)_{n\in \N}$ \\ 
$(u_n)_{n\in\N}\multimap^\omega (v_n)_{n\in\N}=(u_n \multimap v_n)_{n\in\N}$ 
\\ 
$!(u_n)_{n\in\N}=(!u_n)_{n\in\N}$. 
\\ 
$(u_n)_{n\in\N}\leq^\omega (v_n)_{n\in\N}$  when $\forall n\in\N\ u_n\leq v_n$
\end{tabular} 

\end{definition}
IMELL formulae are interpreted as before, but with the operations on $M^\omega$.

Interpretations of an IMELL-theory T are given by morphisms of doctrines: a morphism from $P$ to $P'$ is a pair $\langle F,f\rangle$, where $F:\C\to\C'$ is a functor between the respective base categories preserving finite products, and $f$ is a natural transformation. If $P,P'$ are IMELL $\epsilon$-doctrines, then morphisms preserve the operations $\otimes,\multimap,!$ and the units.
\begin{definition}
    An \emph{interpretation of an IMELL theory T} in $P$ is a morphism $\langle I,i\rangle:\LT_T\to P$.
\end{definition}
Therefore, every IMELL $\epsilon$-doctrine $P_U$ as constructed above provides an interpretation of IMELL + $\epsilon$. By construction, the soundness of this semantics for IMELL follows.

\begin{proposition}
    If an IMELL + $\epsilon$ formula is provable in $T$, then it is satisfied by every interpretation of $T$.
\end{proposition}

%The doctrine 
%substitution in the reverse direction. 

%$C[x,y]\stackrel{t_1(x,y)/w,t_2(x,y)/v,t_3(x,y)/w}\longrightarrow C[u,v,w]$

%The base category B does have products. 

%\sousparagraphe{Posets: ordered monoids with a co-monad}
%An hyperdoctrine is a contravariant functor from B to the category Pos of posets, but the objects of the base category will always be mapped to ordered monoids endowed with a co-monad. 

%We ask for ...

%Injections 

%Projections 

%Exists

%\paragraph{Epsilon}

%The doctrine is said to admit epsilon if the base category B  is endowed with a morphism from to for any formula $\phi[x,y]$
%in such a way that the diagram below commutes: 

\paragraph{Conclusion}
An intriguing question concerns $\tau$, the subnector such that 
$\phi[\tau_x \phi]$ expresses $\forall x\,\phi[x]$. In a classical setting, both operators are not required, since $\tau_x \phi(x)=\epsilon_x \lnot\phi(x)$. The same identification is possible in certain intermediate logics, namely intuitionistic logic enriched with additional axioms.  How, then, could we define $\tau$, the dual of $\epsilon$, within our intuitionistic linear logic equipped with $\epsilon$? In intuitionistic doctrines, universal quantification $\forall$ is available as an adjunction to projection, just as existential quantification $\exists$ is.  So we can define $\tau$ as we defined $\epsilon$. 

Why, then, did we not introduce $\tau$ explicitly? The reason is that, with the ordered monoid we have chosen, we can accommodate either $\epsilon$ or $\tau$, \emph{but not both}. Indeed, every supremum is a maximum, that is, it is attained by at least one value. This guarantees the existence of the existential adjoint. However, it is not the case that every infimum is a minimum. For example, for the set
$\{1/{2^n}\mid n\in\N\}$, the infimum is zero, which does not belong to the set.  

To support both operators, one would need an ordered monoid equipped with a comonad satisfying this property in both directions. A finite ordered monoid would be a possible candidate, but to our knowledge no structure with all the required properties is currently known.  

Since we do have $\bot$, we can define $\phi^\perp$, the linear negation of $\phi$, as $\phi\multimap\bot$, and then define $\tau_x \phi(x)$ by $\epsilon_x \lnot\phi(x)=\epsilon_x (\phi(x)\multimap\bot)$. This possibility deserves further investigation. However, unless $\chi$ is $\bot$, $\chi^\perp$ always collapses to $\bot$, so this approach does not yield a sufficiently fine-grained interpretation.  

This issue is therefore left for further research.

\vspace*{-3,6ex} 

\label{sect:bib}
\bibliographystyle{plain}

\footnotesize
\bibliography{main}

@misc
{
	easychair,
	author       = {Andrei Voronkov},
	title        = {{EasyChair} Conference System},
	year         = {2004--2011},
	howpublished = {\url{easychair.org}}
}

@article{Pasquali2016EpsilonIntuitionistic,  author    = {Fabio Pasquali},  title     = {A Categorical Interpretation of the Intuitionistic, Typed, First Order Logic with {H}ilbert’s $\varepsilon$-Terms},  journal   = {Logica Universalis},  volume    = {10},  number    = {4},  pages     = {407--418},  year      = {2016},  publisher = {Springer},  doi       = {10.1007/s11787-016-0155-y}}

@article{Pasquali2017EpsilonDoctrines,  author  = {Fabio Pasquali},  title   = {{H}ilbert’s Epsilon-Operator in Doctrines},  journal = {IfCoLog Journal of Logics and their Applications},  volume  = {4},  number  = {2},  pages   = {383--402},  year    = {2017},  publisher = {College Publications}}

@article{ChatzikyriakidisPasqualiRetore2017Introduction,  author    = {Stergios Chatzikyriakidis and Fabio Pasquali and Christian Retor{\'e}},  title     = {Introduction to {``{H}ilbert's Epsilon and Tau in Logic, Informatics and Linguistics''}},  journal   = {IfCoLog Journal of Logics and their Applications},  volume    = {4},  number    = {2},  pages     = {221--230},  year      = {2017},  publisher = {College Publications}}

@article{DagninoPasquali2025Lipschitz,  author  = {Francesco Dagnino and Fabio Pasquali},  title   = {Quantitative Equality in Substructural Logic via {L}ipschitz Doctrines},  journal = {Logical Methods in Computer Science},  volume  = {21},  number  = {1},  year    = {2025}}

@article{FleuryQuatrini2004FirstOrderLudics,  author  = {Marie{-}Ren{\'e}e Fleury and Myriam Quatrini},  title   = {First Order in Ludics},  journal = {Mathematical Structures in Computer Science},  volume  = {14},  number  = {2},  pages   = {189--213},  year    = {2004},  doi     = {10.1017/S0960129503004134}}

@phdthesis{Manara2025LinearLogicPiCalculus,  author  = {Giulia Manara},  title   = {Linear Logic: Parallel Cut Elimination and Computation-as-Deduction for the Pi-Calculus},  school  = {Universit{\'e} Paris Cit{\'e} and Universit{\`a} degli Studi Roma Tre},  year    = {2025},  address = {Paris and Rome},  url     = {https://theses.hal.science/tel-05449051}}

@inproceedings{GonthierAbadiLevy1992OptimalReduction,  author    = {Georges Gonthier and Mart{\'\i}n Abadi and Jean{-}Jacques L{\'e}vy},  title     = {The Geometry of Optimal Lambda Reduction},  booktitle = {Proceedings of the 19th ACM SIGPLAN-SIGACT Symposium on Principles of Programming Languages (POPL '92)},  publisher = {ACM},  pages     = {15--26},  year      = {1992},  doi       = {10.1145/143165.143169}}

@incollection{Pitts1992CategoricalLogic,  author    = {Andrew M. Pitts},  title     = {Categorical Logic},  booktitle = {Handbook of Logic in Computer Science},  volume    = {1},  editor    = {Samson Abramsky and Dov M. Gabbay and Thomas S. E. Maibaum},  publisher = {Oxford University Press},  year      = {1992},  pages     = {39--128}}

@article{MaiettiRP,
author = {Maria Emilia Maietti and Fabio Pasquali and Giuseppe Rosolini},
title = {Triposes, exact completions, and {H}ilbert's epsilon-operator},
volume = {10},
journal = {Tbilisi Mathematical Journal},
number = {3},
publisher = {Tbilisi Centre for Mathematical Sciences},
pages = {141 -- 166},
year = {2017},
doi = {10.1515/tmj-2017-0106},
URL = {doi.org/10.1515/tmj-2017-0106}
}

@book{Girard2011blindspot,
	author = {Jean-Yves Girard},
	publisher = {European Mathematical Society},
	title = {The blind spot -- lectures on logic},
	year = {2011}}

@article{Girard2016TranscendentalSyntaxI,  author  = {Jean{-}Yves Girard},  title   = {Transcendental Syntax {I}: Deterministic Case},  journal = {Mathematical Structures in Computer Science},  year    = {2016},  doi     = {10.1017/S0960129515000407}}

@book{HBvol2,
	author = {Hilbert, David and Bernays, Paul},
	booktitle = {Grundlagen der {M}athematik. Bd. 2.},
	language = {German},
	note = {Traduction fran{\c c}aise de F. Gaillard, E. Guillaume et M. Guillaume, L'Harmattan, 2001.},
	publisher = {Springer},
	title = {Grundlagen der {M}athematik. Bd. 2.},
	year = {1939}}

@article{Lawvere1969Adjointness,  author  = {William F. Lawvere},  title   = {Adjointness in Foundations},  journal = {Dialectica},  volume  = {23},  pages   = {281--296},  year    = {1969}}

@article{Gir87,
	author = {Jean-Yves Girard},
	journal = {Theoretical Computer Science},
	number = {1},
	page = {1-102},
	pages = {1--102},
	title = {Linear Logic},
	volume = {50},
	year = {1987}}

@article{heijltjes_proof_2019,
        title = {Proof {Nets} for {First}-{Order} {Additive} {Linear} {Logic}},
        volume = {131},
        issn = {1868-8969},
        doi = {10.4230/LIPICS.FSCD.2019.22},
        language = {en},
        urldate = {2026-06-30},
        journal = {LIPIcs, Volume 131, FSCD 2019},
        author = {Heijltjes, Willem B. and Hughes, Dominic J. D. and Straßburger, Lutz},
        editor = {Geuvers, Herman},
        year = {2019},
        pages = {22:1--22:22}
}

%------------------------------------------------------------------------------
%\appendix

%------------------------------------------------------------------------------
% Index
%\printindex

%------------------------------------------------------------------------------
\end{document}